\documentclass[12pt]{article}
\usepackage{amsmath}
\usepackage{graphicx}
\usepackage{CJK}
\usepackage{multirow}
\usepackage{makecell}

\usepackage{subfigure}
\usepackage[justification=centering]{caption}
\input{amssym.tex}

\def\bc{\begin{center}}       \def\ec{\end{center}}
\def\ba{\begin{array}}        \def\ea{\end{array}}
\def\be{\begin{equation}}     \def\ee{\end{equation}}
\def\bea{\begin{eqnarray}}    \def\eea{\end{eqnarray}}
\def\beaa{\begin{eqnarray*}}  \def\eeaa{\end{eqnarray*}}

\def\r{\right}
\def\lf{\left}

\def\mathbb{\Bbb}

\begin{document}
\baselineskip 18pt
\centerline{\Large\bf Bifurcation of limit cycles from the global center of a class of }
\vskip 0.2 true cm
\centerline{\Large\bf integrable non-Hamilton system under perturbations }
\vskip 0.2 true cm
\centerline{\Large\bf of piecewise smooth polynomials }
\vskip 0.3 true cm

\centerline{\bf  Shiyou Sui, Liqin Zhao$^{*}$}
 \centerline{ School of Mathematical
Sciences, Beijing Normal University,} \centerline{Laboratory of Mathematics and Complex Systems, Ministry of
Education,} \centerline{Beijing 100875, The People's Republic of China}

\footnotetext[1]{
* Corresponding author.~~~This work was supported by NSFC(11671040).\\
E-mail:  zhaoliqin@bnu.edu.cn (L. Zhao), suisy@mail.bun.edu.cn(S. Sui).}

\vskip 0.2 true cm
\noindent{\bf Abstract} In this paper, we perturb the global center of the planar polynomial vector fields $\mathcal{X}(x,y)=(-y(x^2+a^2),x(x^2+a^2))$ ($a\neq0$) inside cubic piecewise smooth polynomials with switching line $y=0$. By using average function of first order, we prove that the sharper bound of the number of limit cycles bifurcating from the period annulus is 6.


\vskip 0.4 true cm
\centerline{{\bf \large $\S 1$. Introduction and the main results}}
\vskip 0.2 true cm

One of the main problems inside the qualitative theory in the qualitative theory  of real planar differential systems is to determine the number of limit cycles which is related to the Hilbert 16th problem [6,11]. A limit cycle is an isolated periodic orbit defined by Poincar\'{e} [8]. There are many phenomena in real world related with the existence of limit cycles, some examples are the van der Pol oscillator [20,21], or the Belousov-Zhabotinskii reaction [1,24]. For more about limit cycles, one can see [2,23].

The notion of a center of real planar differential system is an isolated equilibrium point having a neighborhood such that all the orbits of this neighborhood are periodic with the unique exception of the equilibrium point, which is defined by Poincar\'{e} [18]. Late on a classic way to obtain limit cycles is perturbing the periodic orbits of a center.

In 1999, Iliev [10] considered the polynomial vector fields
$\mathcal{X}(x,y)=(-y,x)$
by perturbing it with polynomials of degree $n$. He studied how many limit cycles can bifurcate from the periodic orbits of the linear center. Later on many people studied the perturbations of polynomial vector fields of the form $\mathcal{X}(x,y)=(-yf(x,y),xf(x,y))$, where $f(0,0)\neq0$, see [13,19,22] and the references they cited.

There are many studies of the limit cycles of continuous and discontinuous piecewise differential systems in $\mathbb{R}^2$ with two pieces separated by a straight line. In general these differential systems are linear, see for instance [5,8,9,16]. T. Carvalho, J. Llibre, and D. J. Tonon [4] studied the number of limit cycles which can bifurcate from the nonlinear center $\mathcal{X}(x,y)=(-y((x^2+y^2)/2)^m,x((x^2+y^2)/2)^m)$, when it is perturbed inside a class of discontinuous piecewise polynomial differential systems of degree $n$ with $k$ pieces. S. Li and Ch. Liu [12] studied the piecewise smooth differential system
$$\mathcal{X}(x,y)=\left\{
\begin{array}{lc}
(-y(1+ax)+\varepsilon P^+(x,y),x(1+ax)+\varepsilon Q^+(x,y)),~~{\text {if}}~x>0,\\
(-y(1+bx)+\varepsilon P^-(x,y),x(1+bx)+\varepsilon Q^-(x,y),~~{\text {if}}~x<0,
\end{array}
\right.$$
where $P^\pm(x,y),~Q^\pm(x,y)$ are polynomials of degree $n$.

 In this paper, we will study the number of limit cycles which can bifurcate from the center of
$$\mathcal{X}(x,y)=(-y(x^2+a^2),x(x^2+a^2)),\eqno{(1.1)}$$
when it is perturbed inside  discontinuous piecewise cubic polynomials, where $a\neq0$. Let $x_1=\displaystyle\frac{1}{a}x,~y_1=\displaystyle\frac{1}{a}y,~t_1=a^2t$, then system (1.1) is transformed into
$$\mathcal{X}(x,y)=(-y(x^2+1),x(x^2+1)),\eqno{(1.2)}$$
here we omit the subscript 1.  Hence, we consider the following perturbations of system (1.2)
$$\mathcal{X}_{\varepsilon}(x,y)=\mathcal{X}(x,y)+\varepsilon\sum\limits_{k=1}^2\mathcal{X}_{S_k}(x,y)(P_k(x,y),Q_k(x,y)),\eqno{(1.3)}$$
where $P_k(x,y)=\sum\limits_{i+j=0}^3a_{ij}^kx^iy^j,~Q_k(x,y)=\sum\limits_{i+j=0}^3b_{ij}^kx^iy^j$, $a_{ij}^k,~b_{ij}^k$ are arbitrary constants,  the characteristic function $\mathcal{X}_S$ of a set $S\subset\mathbb{R}^2$ is defined by
$$\mathcal{X}_S(x,y)=\left\{\begin{array}{lc}
1,~~{\rm if} ~(x,y)\in S,\\
0,~~{\rm if} ~(x,y)\notin S,
\end{array}
\r.$$
and $S_1=\{(x,y)\left|y>0\right.\},~~S_2=\{(x,y)\left|y<0\right.\}$.

The system $(1.3)_{\varepsilon=0}$ has a periodic annulus surround the origin. Then, using average function of first order (see Section 2), we find the maximum number of limit cycles of system (1.3). The main results of this paper is the following.

\noindent{\bf Theorem 1.1.} Suppose that the average function of first order associated to the discontinuous piecewise polynomial differential system (1.3) is non-zero. Then for $\left|\varepsilon\right|>0$ sufficiently small the sharper bound of the number of limit cycles of system (1.3) is  6.

\noindent{\bf Corollary 1.2.} Under the assumption of theorem 1.1, if $a_{ij}^1=a_{ij}^2,~b_{ij}^1=b_{ij}^2$, for $\left|\varepsilon\right|>0$ sufficiently small the sharper bound of the number of limit cycles of system (1.3) is 3.

\noindent{\bf Remark 1.3.} Using the results of [22], we know that system $(1.3)_{\varepsilon=0}$ under the perturbation of smoth polynomials with degree $n$ has at most $2n+2-(-1)^n$.

This paper is organized as follows. In Section 2, we introduce averaging theory for computing periodic solution and extend complete Chebyshev system for studying the number of zeros of average function. The main results is proved in Section 3.

\vskip 0.2 true cm
\centerline{{\bf \large $\S 2$. Preliminary results}}
\vskip 0.2 true cm
In this section we summarize the main tools that we will use to study the bifurcation of limit cycles for system (1.3). First, we introduce the averaging theory for discontinuous piecewise differential systems. The following results stated on the averaging theory are valid for discontinuous piecewise polynomial vector field defined in $\mathbb{R}^n$ and are proved in [15], but we shall state them for our discontinuous piecewise polynomial vector field (1.3) in polar coordinates.

Consider a non-autonomous discontinuous piecewise vector field
$$\frac{{\rm d}r}{{\rm d}\theta}=\mathcal{X}(\theta,r)=\varepsilon F(\theta,r)+\varepsilon^2R(\theta,r,\varepsilon),\eqno{(2.1)}$$
where $r\in\mathbb{R}$, $\theta\in\mathbb{R}/(2\pi\mathbb{Z})$ and
$$F(\theta,r)=\sum\limits_{i=1}^k\mathcal{X}_{S_i}(\theta)F_i(\theta,r),~~~R(\theta,r,\varepsilon)=\sum\limits_{i=1}^k\mathcal{X}_i(\theta)R_i(\theta,r,\varepsilon),$$
where $F_i:\mathbb{S}^1\times D\rightarrow \mathbb{R}^2$, $R_i:\mathbb{S}^1\times D\times (-\varepsilon_0,\varepsilon_0)\rightarrow \mathbb{R}^2$ for $i=1,\ldots,k$ are continuous functions, $2\pi$-periodic in the variable $\theta$, and $D$ is an open interval of $\mathbb{R}$. Here the $\mathbb{S}_i$ are the open intervals $(\theta_i,\theta_{i+1})$ for $i=1,\ldots,k$ and $0\leq\theta_1\leq\cdots\leq\theta_k\leq2\pi\leq\theta_{k+1}=\theta_1+2\pi$. We define
$$D_rF(\theta,r)=\sum\limits_{i=1}^k\mathcal{X}_{S_i}(\theta,r)D_rF_i(\theta,r).$$
The {\it average function} $f:D\rightarrow \mathbb{R}$ is defined by
$$f(r)=\int_0^TF(\theta,r){\rm d}\theta.$$
We recall that if $r(\theta,r_0)$ is the solution of the vector field $\mathcal{X}(\theta,r)$ such that $r(0,r_0)=r_0$, then we have
$$r(2\pi,r_0)-r_0=\varepsilon f(r)+O(\varepsilon^2).$$
So for $\varepsilon>0$ suffieiently small the simple zeros of the averaged function $f(r)$ provides limit cycles of the vector field $\mathcal{X}(\theta,r)$.

In the next result we present a version of the averaging theory for discontinuous piecewise vector fields, that is proved in [15], adapted to  differential equation (2.1).  We note that in [15] the averaging theory uses that the Brouwer degree of a function $f$ in a neighborhood of a zero $\bar{r}$ of the function $f(r)$ is non-zero, while here we substitute this condition saying that the zero $\bar{r}$ is simple (i.e. $\frac{{\rm d}f}{{\rm d}r}(\bar{r})\neq0$), because this last condition implies that the mentioned Brouwer degree  non-zero. See for more details [3,17].

\noindent{\bf Lemma 2.1.} Assume that the following conditions hold for the discontinuous piecewise vector field $\mathcal{X}(\theta,r)$.\\
(i) For $i=1,\ldots,k$ the functions $F_i(\theta,r)$ and $R_i(\theta,r)$ are locally Lipschitz with respect to $r$, and $2\pi$-periodic with respect to $\theta$.\\
(ii) Let $\bar{r}\in D$ be a simple zero of the average function $f(r)$.\\
Then for $\varepsilon>0$ sufficiently small, there exists a $2\pi$-periodic solution $r(\theta,\varepsilon)$ of the vector field $\mathcal{X}(\theta,r)$ such that $r(0,\varepsilon)\rightarrow \bar{r}$ as $\varepsilon\rightarrow 0$.

In order to study the number of zeros of the averaging function we will use the following results.

\noindent{\bf Definition 2.2.}([14]) Let $\mathbb{U}$ be a set and let $f_1,f_2,\ldots,f_n:\mathbb{U}\rightarrow\mathbb{R}$, we say that $f_1,\ldots,f_n$ are linearly independent functions if and only if we have that
$$\sum\limits_{i=1}^nk_if_i(x)=0~~{\text{for all}}~x\in \mathbb{U}~~\Rightarrow~~~k_1=k_2=\ldots=k_n=0.$$

\noindent{\bf Lemma 2.3.}([14]) If $f_1,f_2,\ldots,f_n:\mathbb{U}\rightarrow\mathbb{R}$ are linearly independent then there exit $x_1,x_2,\ldots,x_{n-1}\in \mathbb{U}$ and $k_1,k_2,\ldots,k_n\in \mathbb{R}$ such that for every $i\in\{1,\ldots,n-1\}$
$$\sum\limits_{j=1}^nk_jf_j(x_i)=0.$$

\noindent{\bf Definition 2.4.}([7]) Let $f_0,f_1,\ldots,f_{n-1}$ be analytic functions on an open interval $L$ of $\mathbb{R}$.\\
(a) $(f_0,f_1,\ldots,f_{n-1})$ is a \textit{Chebyshev system} (in short, T-system) on $L$ if any nontrivial linear combination
$$\alpha_0f_0(x)+\alpha_1f_1(x)+\ldots+\alpha_{n-1}f_{n-1}(x)$$
has at most $n-1$ isolated zeros on $L$.\\
(b) $(f_0,f_1,\ldots,f_{n-1})$ is an \textit{complete Chebyshev system} (in short, CT-system) on $L$ if $(f_0,f_1,\ldots,f_{k-1})$ is a T-system for all $k=1,2,\ldots,n$.\\
(c) $(f_0,f_1,\ldots,f_{n-1})$ is an \textit{extend complete Chebyshev system} (in short, ECT-system) on $L$ if, for all $k=1,2,\ldots,n$, any nontrivial linear combination
$$\alpha_0f_0(x)+\alpha_1f_1(x)+\ldots+\alpha_{k-1}f_{k-1}(x)$$
has at most $k-1$ isolated zeros on $L$ counted with multiplicities.

It is clear that if $(f_0,f_1,\ldots,f_{n-1})$ is an EXT-system on $L$, then $(f_0,f_1,\ldots,f_{n-1})$ is a CT-system on $L$. However, the reverse implication is not true.

\noindent{\bf Lemma 2.5.}([7]) $(f_0,f_1,\ldots,f_{n-1})$ is an ECT-system on $L$ if and only if for each $k=1,2,\ldots,n$ the continuous Wronskian of $(f_0,f_1,\ldots,f_{k-1})$ at $x\in L$ is not zero, that is
$$W[f_0,f_1,\dots,f_{k-1}](x)=\left|\begin{matrix}
f_0(x) & f_1(x) & \cdots & f_{k-1}(x)\\
f_0'(x)&f_1'(x)&\cdots&f_{k-1}'(x)\\
\cdots&\cdots&\cdots&\cdots\\
f_0^{(k-1)}(x)&f_1^{(k-1)}(x)&\cdots&f_{k-1}^{(k-1)}(x)
\end{matrix}\right|\neq 0~~(x\in L).
$$

By Definition 2.2 and Lemma 2.5, it is easy to get the following.

\noindent{\bf Proposition 2.6.} If $(f_0,f_1,\ldots,f_{n-1})$ is an ECT-system on $L$, then they are linearly independent.

\vskip 0.4 true cm
\centerline{{\bf \large $\S 3$. Proof of main results}}
\vskip 0.2 true cm

Using the change to polar coordinates $x=r\cos\theta,~y=r\sin\theta$, we transform the differential system (1.3) into
$${\begin{split}
{\dot{r}}&=\varepsilon\sum\limits_{k=1}^2\left[\cos\theta(\mathcal{X}_{S_k}\cdot P_k)(r\cos\theta,r\sin\theta)+\sin\theta(\mathcal{X}_{S_k}\cdot Q_k)(r\cos\theta,r\sin\theta)\right],\\
{\dot{\theta}}&=r^2\cos^2\theta+\frac{\varepsilon}{r}\sum\limits_{k=1}^2\left[\cos\theta(\mathcal{X}_{S_k}\cdot Q_k)(r\cos\theta,r\sin\theta)-\sin\theta(\mathcal{X}_{S_k}\cdot P_k)(r\cos\theta,r\sin\theta)\right].
\end{split}}$$
Taking $\theta$ as the new independent variable the previous differential system becomes the differential equation
$$\displaystyle\frac{{\rm d}r}{{\rm d}\theta}=\varepsilon F(\theta,r)+\varepsilon^2R(\theta,r,\varepsilon),\eqno{(3.1)}$$
where
$$F(\theta,r)=\sum\limits_{k=1}^2\frac{\cos\theta (\mathcal{X}_{S_k}\cdot P_k)(r\cos\theta,r\sin\theta)+\sin\theta (\mathcal{X}_{S_k}\cdot Q_k)(r\cos\theta,r\sin\theta)}
{r^2\cos^2\theta+1},$$
$${\begin{split}
R(\theta,r,\varepsilon)&=-\sum\limits_{k=1}^2\frac{\cos\theta (\mathcal{X}_{S_k}\cdot P_k)(r\cos\theta,r\sin\theta)+\sin\theta (\mathcal{X}_{S_k}\cdot Q_k)(r\cos\theta,r\sin\theta)}
{r^2\cos^2\theta+1}\\
&\cdot\frac{\cos\theta(\mathcal{X}_{S_k}\cdot Q_k)(r\cos\theta,r\sin\theta)-\sin\theta(\mathcal{X}_{S_k}\cdot P_k)(r\cos\theta,r\sin\theta)}{r(r^2\cos^2\theta+1)+\varepsilon\left[\cos\theta(\mathcal{X}_{S_k}\cdot Q_k)(r\cos\theta,r\sin\theta)-\sin\theta(\mathcal{X}_{S_k}\cdot P_k)(r\cos\theta,r\sin\theta)\right]}.\end{split}}$$
Therefore, the average function of first order is
$${\begin{split}
f(r)&=\int_0^{2\pi}F(\theta,r){\rm d}\theta\\
&=\sum\limits_{i+j=0}^3\int_0^{\pi}\frac{a_{ij}^1r^{i+j}\cos^{i+1}\theta\sin^j\theta+b_{ij}^1r^{i+j}\cos^{i}\theta\sin^{j+1}\theta}
{r^2\cos^2\theta+1}{\rm d}\theta\\
&~+\sum\limits_{i+j=0}^3\int_{\pi}^{2\pi}\frac{a_{ij}^2r^{i+j}\cos^{i+1}\theta\sin^j\theta+b_{ij}^2r^{i+j}\cos^{i}\theta\sin^{j+1}\theta}
{r^2\cos^2\theta+1}{\rm d}\theta\\
&=\sum\limits_{i=0}^3r^i\sum\limits_{j=0}^{i+1}\omega_{i+1-j,j}^1\int_0^{\pi}\frac{\cos^{i+1-j}\theta\sin^j\theta}
{r^2\cos^2\theta+1}{\rm d}\theta\\
&~+\sum\limits_{i=0}^3r^i\sum\limits_{j=0}^{i+1}\omega_{i+1-j,j}^2\int_{\pi}^{2\pi}\frac{\cos^{i+1-j}\theta\sin^j\theta}
{r^2\cos^2\theta+1}{\rm d}\theta
\end{split}}$$
where $\omega_{i,j}^k=a_{i-1,j}^k+b_{i,j-1}^k,~1\leq i+j\leq 4,~a_{-1,j}^k=b_{i,-1}^k=0,~k=1,2$. Note that $\omega_{i,j}^k$ are also arbitrary, since $a_{i,j}^k,~b_{i,j}^k$ are arbitrary.

In order to simplify the notation, we define the following functions:
$$I_{i,j}(r)=\int_0^{\pi}\frac{\cos^i\theta\sin^j\theta}{r^2\cos^2\theta+1}{\rm d}\theta,$$
$$J_{i,j}(r)=\int_{\pi}^{2\pi}\frac{\cos^i\theta\sin^j\theta}{r^2\cos^2\theta+1}{\rm d}\theta.$$
Then, it is easy to check that
$$J_{i,j}(r)=(-1)^{i+j}I_{i,j}(r).$$
Notice that in the interval $(0,+\infty)$, the zeros of the function $f(r)$ coincide with the zeros of the function $F(r)=rf(r)$. Therefore, in order to simplify further computation, we will study the function $F(r)$ instead of $f(r)$.

Using above notation, we can obtain that
\begin{align}
F(r)&=rf(r)=\sum\limits_{i=1}^4r^i\sum\limits_{j=0}^i\omega_{i-j,j}^1I_{i-j,j}(r)+\sum\limits_{i=1}^4r^i\sum\limits_{j=0}^i\omega_{i-j,j}^2J_{i-j,j}(r)\notag\\
&=\sum\limits_{i=1}^4r^i\sum\limits_{j=0}^i\left(\omega_{i-j,j}^1+(-1)^{i}\omega_{i-j,j}^2\right)I_{i-j,j}(r)\notag\\
&=\sum\limits_{i=1}^4r^i\sum\limits_{j=0}^i\mu_{i-j,j} I_{i-j,j}(r)\tag{3.2}
\end{align}
where $\mu_{i-j,j}=\omega_{i-j,j}^1+(-1)^{i}\omega_{i-j,j}^2$. Note that $\mu_{i-j,j}$ are independent, since $\omega_{i-j,j}^1,~\omega_{i-j,j}^2$ are arbitrary.
By direct computation we have the following:
$$\lf\{\begin{array}{lc}
I_{1,0}(r)=I_{1,1}(r)=I_{3,0}(r)=I_{1,2}(r)=I_{3,1}(r)=I_{1,3}(r)=0,\\
I_{0,0}(r)=\frac{\pi}{\sqrt{1+r^2}},~~I_{0,1}(r)=2\frac{\arctan r}{r},\\
I_{2,0}(r)=\frac{\pi}{r^2}(1-\frac{1}{\sqrt{1+r^2}}),~~I_{2,1}(r)=\frac{2}{r^3}(r-\arctan r),\\
I_{0,2}(r)=I_{0,0}(r)-I_{2,0}(r),~~I_{0,3}(r)=I_{0,1}(r)-I_{2,1}(r),~~I_{4,0}(r)=\frac{\pi}{2r^2}-\frac{1}{r^2}I_{2,0}(r),\\
I_{2,2}(r)=I_{2,0}(r)-I_{4,0}(r),~~I_{0,4}(r)=I_{4,0}(r)-2I_{2,0}(r)+I_{0,0}(r).
\end{array}
\r.\eqno{(3.3)}$$

Substituting (3.3) into (3.2), we have that
\begin{align}
F(r)&=\mu_{0,1}rI_{0,1}(r)+r^2(\mu_{2,0}I_{2,0}(r)+\mu_{0,2}I_{0,2}(r))+r^3(\mu_{2,1}I_{2,1}(r)+\mu_{0,3}I_{0,3}(r))\notag\\
&~~+r^4(\mu_{4,0}I_{4,0}(r)+\mu_{2,2}I_{2,2}(r)+\mu_{0,4}I_{0,4}(r))\notag\\
&=(\mu_{0,2}r^2+\mu_{0,4}r^4)I_{0,0}(r)+(\mu_{0,1}r+\mu_{0,3}r^3)I_{0,1}(r)+((\mu_{2,0}-\mu_{0,2})r^2+(\mu_{2,2}-2\mu_{0,4})r^4)I_{2,0}(r)\notag\\
&~~+(\mu_{2,1}-\mu_{0,3})r^3I_{2,1}(r)+(\mu_{4,0}-\mu_{2,2}+\mu_{0,4})r^4I_{4,0}\notag\\
&=2\nu_1r+\pi\nu_2r^2+\pi \nu_3\frac{r^2}{\sqrt{1+r^2}}+\pi\nu_4\frac{r^4}{\sqrt{1+r^2}}+2\nu_5\arctan r \notag\\
&~~+2\nu_6r^2\arctan r  +\pi\nu_7(1-\frac{1}{\sqrt{1+r^2}})\tag{3.4}
\end{align}
where

$$\lf\{\begin{array}{lc}
\nu_1=\mu_{2,1}-\mu_{0,3},\\
\nu_2=\frac{1}{2}\mu_{4,0}+\frac{1}{2}\mu_{2,2}-\frac{3}{2}\mu_{0,4},\\
\nu_3=\mu_{0,2}-\mu_{2,2}+2\mu_{0,4},\\
\nu_4=\mu_{0,4},\\
\nu_5=\mu_{0,1}-\mu_{2,1}+\mu_{0,3},\\
\nu_6=\mu_{0,3},\\
\nu_7=\mu_{2,0}-\mu_{0,2}-\mu_{4,0}+\mu_{2,2}-\mu_{0,4}.
\end{array}
\r.$$

It follows from direct computation that
$${\rm det}\frac{\partial(\nu_1,\nu_2,\nu_3,\nu_4,\nu_5,\nu_6,\nu_7)}{\partial(\mu_{2,1},\mu_{4,0},\mu_{0,2},\mu_{0,4},\mu_{0,1},\mu_{0,3},\mu_{2,0})}=\frac{1}{2}\neq0.$$
By the arbitrariness of $\mu_{ij}$, we have that $\nu_i$ $(i=1,\ldots,7)$ are independent. Hence, the generating functions of $F(r)$ are the following:
$${\begin{split}
&f_1(r):=r,~~f_2(r):=r^2,~~f_3(r):=\frac{r^2}{\sqrt{1+r^2}},~~f_4(r):=\frac{r^4}{\sqrt{1+r^2}},\\
&f_5(r):=\arctan r,~~f_6(r):=r^2\arctan r,~~f_7(r):=1-\frac{1}{\sqrt{1+r^2}}.
\end{split}}$$

Next, we will prove that $(f_1(r),\ldots,f_7(r))$ is an ECT-system. We introduce the notation
$W_k(r)=W[f_1(r),f_2(r),\ldots,f_k(r)]$. Direct computation, we have
$$W_2(r)=r^2,~~W_3(r)=-\frac{3r^4}{(r^2+1)^{\frac{5}{2}}},~~W_4(r)=-\frac{6r^7(4r^2+5)}{(r^2+1)^5}.$$
It is obvious that $W_k(r)\neq0~(k=2,3,4)$ on $r\in(0,+\infty)$.
\begin{figure*}[!hbt]
\centering
{\includegraphics[scale=0.4]{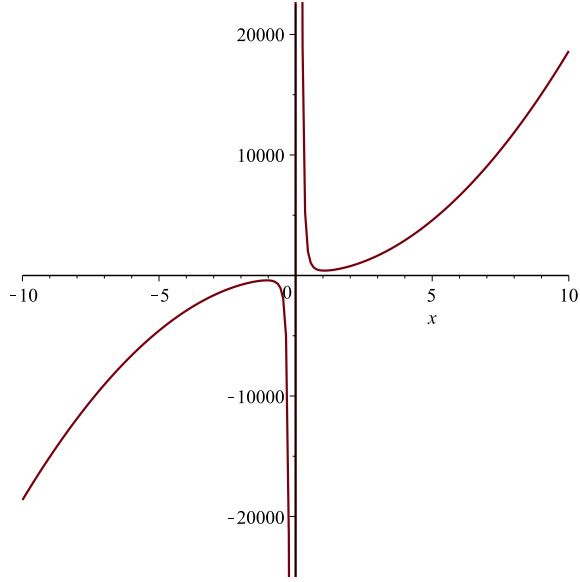}}
\caption{The curve of $\frac{g_1'(r)}{r^7}$.}
\end{figure*}

For $k=5$, using mathematical soft such as Maple, we get that
$$W_5(r)=\frac{12r^3}{(r^2+1)^9}g_1(r),$$
where
$${\begin{split}
g_1(r)&=12r^{10}\arctan r-27r^8\arctan r+12r^9-258r^6\arctan r+169r^7\\
&~-507r^4\arctan r+411r^5-408r^2\arctan r+373r^3-120\arctan r+120r.
\end{split}}$$
Then,
$${\begin{split}
g_1'(r)&=r(120r^8\arctan r-216r^6\arctan r+120r^7-1548r^4\arctan r+1144r^5\\
&-2028r^2\arctan r+1836r^3-816\arctan r+831r).
\end{split}}$$
From the curve (Figure 1) of $\frac{g_1'(r)}{r^7}$, we know that $g_1'(r)>0$ on $(0,+\infty)$.
So, we get that $g_1(r)>g_1(0)=0$. Hence, $W_5(r)\neq0$ on $(0,+\infty)$.

When $k=6$, we have
$$W_6(r)=-\frac{24r}{(r^2+1)^{13}}g_2(r),$$
where
$${\begin{split}
g_2(r)&=216r^{12}\arctan r+168r^{10}\arctan r+216r^{11}+843r^8\arctan r+96r^9\\
&+4986r^6\arctan r-3061r^7+8175r^4\arctan r-6655r^5\\
&+5160r^2\arctan r-4800r^3+1080\arctan r-1080r.
\end{split}}$$
By the curve (see Figure 2) of $\frac{g_2(r)}{r^{10}}$, we have that $W_6(r)\neq0$ on $(0,+\infty)$.
\begin{figure*}[!!ht]
\centering
{\includegraphics[scale=0.4]{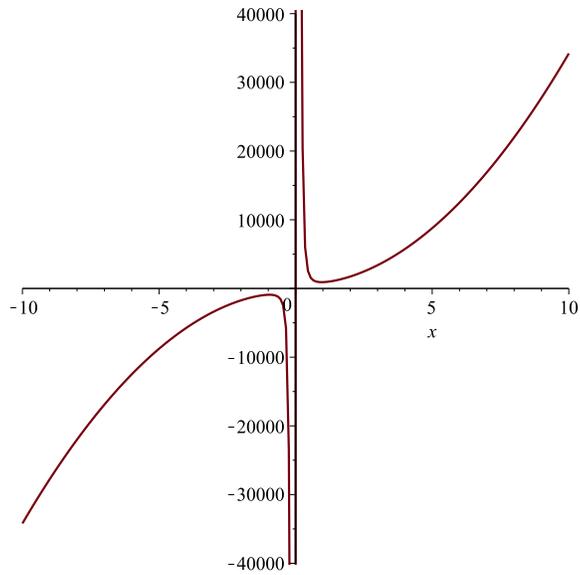}}
\caption{The curve of $\frac{g_2(r)}{r^{10}}$.}
\end{figure*}

When $k=7$, we have
$$W_7(r)=\frac{1728}{(x^2+1)^{\frac{35}{2}}}g_3(r),$$
where
$${\begin{split}
g_3(r)&=3120r^{11}\arctan r-4864r^{10}\sqrt{r^2+1}+10500r^9\arctan r+3120r^{10}-14048r^8\sqrt{x^2+1}\\
&+13155r^7\arctan r+9460r^8-14224r^6\sqrt{r^2+1}+7350r^5\arctan r+10279 r^6-6020r^4\sqrt{x^2+1}\\
&+1575r^3\arctan r+4985r^4-1120r^2\sqrt{r^2+1}+1200r^2-160\sqrt{r^2+1}+160.
\end{split}}$$
\begin{figure*}[!!ht]
\centering
{\includegraphics[scale=0.4]{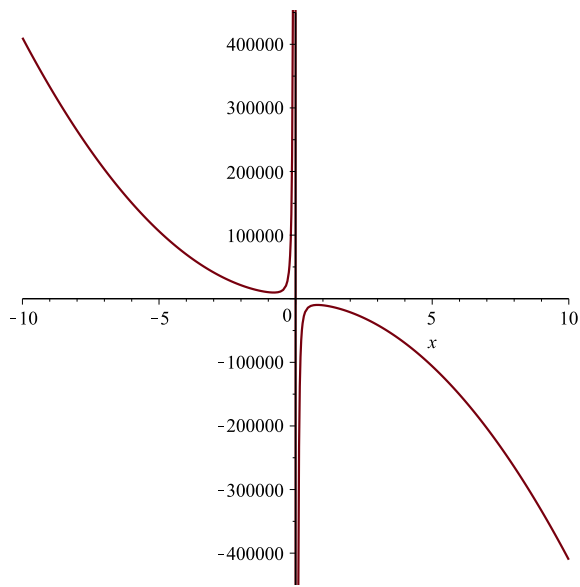}}
\caption{The curve of $\frac{g_{31}(r)}{r^{8}}$.}
\end{figure*}
Then,
$${\begin{split}
g_3'(r)&=\frac{r}{\sqrt{r^2+1}}\left(34320r^9\sqrt{r^2+1}\arctan r-53504r^{10}+94500r^7\sqrt{r^2+1}\arctan r\right.\\
&+34320r^8\sqrt{r^2+1}-175072r^8+92085r^5\sqrt{x^2+1}\arctan r+83060r^6\sqrt{r^2+1}\\
&-211952r^6+36750r^3\sqrt{r^2+1}\arctan r+67449r^4\sqrt{r^2+1}\\
&-115444r^4+4725r\sqrt{r^2+1}\arctan r+21515r^2\sqrt{r^2+1}\\
&\left.-27440r^2+2400\sqrt{r^2+1}-2400\right)
\end{split}}$$
and
$$\frac{{\rm d}}{{\rm d}r}\left(\frac{g_3'(r)\sqrt{r^2+1}}{r}\right)=-\frac{1}{\sqrt{r^2+1}}g_{31}(r)$$
where
$${\begin{split}
g_{31}(r)&=-343200r^{10}\arctan r+535040r^9\sqrt{r^2+1}-1064880r^8\arctan r-343200r^9\\
&+1400576r^7\sqrt{r^2+1}-1214010r^6\arctan r-950480r^7+1271712r^5\sqrt{r^2+1}\\
&-607425r^4\arctan r-927690r^5+461776r^3\sqrt{r^2+1}-119700r^2\arctan r\\
&-371091r^3+54880r\sqrt{r^2+1}-4725\arctan r-50155r.
\end{split}}$$
From the curve (see Figure 3) of $\frac{g_{31}(r)}{r^8}$, we know $g_{31}(r)<0$ on $(0,+\infty)$.

Hence, $\frac{{\rm d}}{{\rm d}r}\left(\frac{g_3'(r)\sqrt{r^2+1}}{r}\right)>0$ on $(0,+\infty)$. Therefore, we have that $\frac{g_3'(r)\sqrt{r^2+1}}{r}>\lim\limits_{r\rightarrow0}\frac{g_3'(r)\sqrt{r^2+1}}{r}=0$. That is equivalent to $g_3'(r)>0$ on $(0,+\infty)$. Thus, $g_3(r)>g_3(0)=0$, which implies $W_7(r)\neq0$ on $(0,+\infty)$.
By above analysis we have proved that $(f_1(r),\ldots,f_7(r))$ is an ECT-system on $(0,+\infty)$. So, by Lemma 2.5, $F(r)$ has at most 6 zeros on $(0,+\infty)$. Using proposition 2.6 and Lemma 2.3, we have that $F(r)$ can have 6 zeros on $(0,+\infty)$. Hence, by Lemma 2.1, theorem 1.1 is proved.

If $a_{ij}^1=a_{ij}^2,~b_{ij}^1=b_{ij}^2$, then $\omega_{i-j,j}^1=\omega_{i-j,j}^2$, which imply $\mu_{i-j,j}=0$ (i is odd). Then, the generating functions of $F(r)$ become
$$f_2(r):=r^2,~~f_3(r):=\frac{r^2}{\sqrt{1+r^2}},~~f_4(r):=\frac{r^4}{\sqrt{1+r^2}},~~
f_7(r):=1-\frac{1}{\sqrt{1+r^2}}.$$
It is easy to check that $(f_2(r),f_3(r),f_4(r),f_7(r))$ is an ECT-system on $(0,+\infty)$. Then, we ends the proof of Corollary 1.2.

\end{document}